\documentclass[a4paper]{amsart}  

\usepackage{amssymb} \usepackage{amscd} \usepackage{times}

\def\zbb{\mathbb{Z}}  
  
  \def\phi{\varphi}
 \def\p1{{\mathbb{P}^1_\zbb}}

\newtheorem{Theorem}{\quad Theorem}[section]

\newtheorem{Proposition}[Theorem]{\quad Proposition}

\newtheorem{Lemma}[Theorem]{\quad Lemma}

\begin{document}

\title{ About Brezis-Merle Problem with Holderian condition: The case of one or two blow-up points.}

\author{Samy Skander Bahoura}

\address{Departement de Mathematiques, Universite Pierre et Marie Curie, 2 place Jussieu, 75005, Paris, France.}
              
\email{samybahoura@yahoo.fr, samybahoura@gmail.com} 

\date{}

\maketitle

\begin{abstract}

We consider the following problem on open set  $ \Omega $ of $ {\mathbb R}^2 $:

\begin{displaymath}  \left \{ \begin {split} 
      -\Delta u_i & = V_i e^{u_i} \,\, &&\text{in} \!\!&&\Omega \subset {\mathbb R}^2, \\
                  u_i  & = 0  \,\,             && \text{in} \!\!&&\partial \Omega.               
\end {split}\right.
\end{displaymath}

We give a quantization analysis of the previous problem under the conditions:

$$ 0 \leq V_i \leq b  < + \infty, \,\, {\rm and}\,\, \int_{\Omega} e^{u_i} dy  \leq C. $$

On the other hand, if  we assume that 

$$ \int_{\Omega} V_i e^{u_i} dy \leq 4\pi, $$

or, $ V_i  $ $ s $-holderian with $ 1/2 < s \leq 1$, and,  %then we have a compactness result, namely:

$$ \int_{\Omega}V_i e^{u_i} dy \leq 24\pi-\epsilon , \,\, \epsilon >0 $$
 
 then we have a compactness result, namely:

$$ \sup_{\Omega} u_i \leq c=c(b, C, A, s, \epsilon, \Omega). $$

where $ A $ is the holderian constant of $ V_i $.

\end{abstract}

\section{Introduction and Main Results} 

We set $ \Delta = \partial_{11} + \partial_{22} $  on open set $ \Omega $ of $ {\mathbb R}^2 $ with a smooth boundary.

\smallskip

We consider the following problem on $ \Omega \subset {\mathbb R}^2 $:

\begin{displaymath} (P) \left \{ \begin {split} 
      -\Delta u_i & = V_i e^{u_i} \,\, &&\text{in} \!\!&&\Omega \subset {\mathbb R}^2, \\
                  u_i  & = 0  \,\,             && \text{in} \!\!&&\partial \Omega.               
\end {split}\right.
\end{displaymath}
We assume that,

$$ \int_{\Omega} e^{u_i} dy  \leq C, $$

and,

$$ 0 \leq V_i \leq b  < + \infty $$

The above equation is called, the Prescribed Scalar Curvature equation, in
relation with conformal change of metrics. The function $ V_i $ is the
prescribed curvature.

Here, we try to find some a priori estimates for sequences of the
previous problem.

Equations of this type were studied by many authors, see [1-26]. One can
see in [8] different results for the solutions of those type of
equations with or without boundaries conditions and, with minimal
conditions on $ V $, for example we suppose $ V_i \geq 0 $ and  $ V_i
\in L^p(\Omega) $ or $ V_ie^{u_i} \in L^p(\Omega) $ with $ p \in [1,
+\infty] $. 

Among other results, we  can see in [8], the following important Theorem,

{\bf Theorem A}{\it (Brezis-Merle [8])}.{\it If $ (u_i)_i $ and $ (V_i)_i $ are two sequences of functions relatively to the previous problem $ (P) $ with, $ 0 < a \leq V_i \leq b < + \infty $, then, for all compact set $ K $ of $ \Omega $,

$$ \sup_K u_i \leq c = c(a, b, K, \Omega). $$}

We can find in [8] an interior estimate if we assume $ a=0 $, but we need an assumption on the integral of $ e^{u_i} $. We have in [8]:

\smallskip

{\bf Theorem B} {\it (Brezis-Merle [8])}.{\it If $ (u_i)_i $ and $ (V_i)_i $ are two sequences of functions relatively to the previous problem $ (P) $ with, $ 0 \leq V_i \leq b < + \infty $, and,

$$ \int_{\Omega} e^{u_i} dy  \leq C, $$

then, for all compact set $ K $ of $ \Omega $,

$$ \sup_K u_i \leq c = c(b, C, K, \Omega). $$}

If, we assume $ V $ with more regularity, we can have another type of estimates, $ \sup + \inf $. It was proved, by Shafrir, see [23], that, if $ (u_i)_i, (V_i)_i $ are two sequences of functions solutions of the previous equation without assumption on the boundary and, $ 0 < a \leq V_i \leq b < + \infty $, then we have the following interior estimate:

$$ C\left (\dfrac{a}{b} \right ) \sup_K u_i + \inf_{\Omega} u_i \leq c=c(a, b, K, \Omega). $$

One can see in [12], an explicit value of $ C\left (\dfrac{a}{b}\right ) =\sqrt {\dfrac{a}{b}} $. In his proof, Shafrir has used a blow-up function, the Stokes formula and an isoperimetric inequality, see [6]. For Chen-Lin, they have used the blow-up analysis combined with some geometric type inequality for the integral curvature.

\bigskip

Now, if we suppose $ (V_i)_i $ uniformly Lipschitzian with $ A $ the
Lipschitz constant, then, $ C(a/b)=1 $ and $ c=c(a, b, A, K, \Omega)
$, see Brezis-Li-Shafrir [7]. This result was extended for
H\"olderian sequences $ (V_i)_i $ by Chen-Lin, see  [12]. Also, we have in [18] an extension of the Brezis-Li-Shafrir result to compact
Riemann surface without boundary. We have in [19] explicit form,
($ 8 \pi m, m\in {\mathbb N}^* $ exactly), for the numbers in front of
the Dirac masses, when the solutions blow-up. Here, the notion of isolated blow-up point is used. Also, we have in [13]  and [26] refined estimates near the isolated blow-up points and the  bubbling behavior  of the blow-up sequences.

\bigskip

In [8], Brezis and Merle proposed the following Problem:

\bigskip

{\bf Problem} {\it (Brezis-Merle [8])}.{\it If $ (u_i)_i $ and $ (V_i)_i $ are two sequences of functions relatively to the previous problem $ (P) $ with, 

$$ 0 \leq V_i \to V \,\, {\rm in } \,\,  C^0(\bar \Omega). $$

$$ \int_{\Omega} e^{u_i} dy  \leq C, $$

Is it possible to prove that:

$$ \sup_{\Omega} u_i \leq c = c(C, V, \Omega) \,\, ? $$}

Here, we assume more regularity on $ V_i $, we suppose that $ V_i \geq 0 $ is $ C^s $ ($ s$-holderian) $ 1/2 < s \leq 1)$ . We give the answer where $ b C < 24 \pi $.

On other hand, in our work we give a complete caracterization of the blow-up analysis on the boundary.

\bigskip

In the similar way, we have in dimension $ n \geq 3 $, with different methods, some a priori estimates of the type $ \sup \times \inf $ for equation of the type:

$$ -\Delta u+ \dfrac{n-2}{4(n-1)}  R_g(x)u = V(x) u^{(n+2)/(n-2)} \,\,\, {\rm on } \,\, M. $$

where $ R_g $ is the scalar curvature of a riemannian manifold $ M $, and $ V $ is a function. The operator $  \Delta=\nabla^i(\nabla_i) $ is the Laplace-Beltrami operator on $ M $.

\smallskip

When $ V \equiv 1 $ and $ M $ compact, the previous equation is the Yamabe equation. T. Aubin and R. Scheon solved the Yamabe problem, see for example [1].  Also, we can have an idea on the Yamabe Problem in [16].  If $ V $ is not a constant function, the previous equation is called a prescribing curvature equation, we have many existence results see also [1]. 

\smallskip

Now, if we look at the problem of a priori bound for the previous equation, we can see in [2, 3, 4, 5, 11, 17, 21] some results concerning the $ \sup \times \inf $ type of inequalities when the manifold $ M $ is  the sphere or more generality a locally conformally flat manifold. For these results, the moving-plane was used, we refer to [9, 14, 20] to have an idea on this method and some applications of this method. 

Also, there are similar problems defined on complex manifolds for the Complex Monge-Ampere equation, see [24, 25], with various inequalitites of type $ \sup+\inf $.

\smallskip

Our main results are:

\begin{Theorem} Assume $ \Omega=B_1(0) $, and,

$$ u_i(x_i)=\sup_{B_1(0)} u_i \to + \infty. $$

 There is a finite number of sequences $ (x_i^k)_i, (\delta_i^k) , 0 \leq k \leq m $, such that:

$$  (x_i^0)_i \equiv (x_i)_i, \,\, \delta_i^0=\delta_i=d(x_i, \partial B_1(0)) \to 0, $$

%$$ (x_i^1)_i \equiv (t_i)_i, \,\, \delta_i^1=\delta'_i=d(t_i, \partial (B_1(0)-B(x_i, \delta_i \epsilon)), $$

and each $\delta_i^k $ is of order $ d(x_i^k, \partial B_1(0)) $.

and,

$$ u_i(x_i^k)=\sup_{B_1(0)- \cup_{j=0}^{k-1}  B(x_i^j, \delta_i^j \epsilon)} u_i \to + \infty, $$

$$ u_i(x_i^k) + 2\log \delta_i^k \to + \infty, $$

$$ \forall \,\, \epsilon >0, \,\, \sup_{B_1(0)- \cup_{j=0}^{m}  B(x_i^j, \delta_i^j \epsilon)} u_i \leq C_{\epsilon} $$

$$ \forall \,\, \epsilon >0, \,\,\, \limsup_{i \to + \infty} \int_{B(x_i^k, \delta_i^k\epsilon)} V_i e^{u_i}  dy \geq 4 \pi >0 .$$

 If we assume:

$$ V_i \to V\,\, {\rm in } \,\, C^0(\bar B_1(0)), $$

then,

$$  \forall \,\, \epsilon >0, \,\,\, \limsup_{i \to + \infty} \int_{B(x_i^k, \delta_i^k\epsilon)} V_i e^{u_i}  dy = 8 \pi m_k, \,\, m_k \in  {\mathbb N}^*.$$

And, thus, we have the following convergence in the sense of distributions:

$$ \int_{B_1(0)} V_i e^{u_i}  dy  \to \int_{B_1(0)} V e^{u}  dy +\sum_{k=0}^m 8 \pi m_k' \delta_{x_0^k} ,\, m_k' \in  {\mathbb N}^*, \,\, x_0^k \in \partial B_1(0). $$

\end{Theorem}

\begin{Theorem}Assume that:

$$ \int_{B_1(0)} V_i e^{u_i}  dy  \leq 4\pi, $$

Then,

$$ u_i(x_i)=\sup_{B_1(0)} u_i \leq c=c(b, C), $$

\end{Theorem}

\begin{Theorem}Assume that, $ V_i  $ is uniformly $ s- $holderian with $ 1/2 < s \leq 1$, and,

$$ \int_{B_1(0)} V_i e^{u_i}  dy  \leq 24\pi-\epsilon, \,\, \epsilon>0,$$

then we have:

$$ \sup_{\Omega} u_i \leq c=c(b, C, A, s, \epsilon, \Omega). $$

where $ A $ is the holderian constant of $ V_i $.

\end{Theorem}

{\bf Question 1:} ({\it a Bartolucci type result; one holderian singularity}): with the same technique, assume that:

$$ V_i(x)=(1+x_1^s) W_i(x) \,\, \text{for \, example \, and } 0 \in \partial \Omega \,\,  $$

with $ W_i $ uniformly lipschitzian and , $ 0 <s \leq 1$, can one conclude with the Pohozaev identity that  the sequence is compact ? here we extend the case $ 0 < s \leq 1/2 $. 

{\bf Question 2: } (the limit case $ s=1/2 $) assume that $ V_i $ is uniformly $ 1/2 $- holderian with $ A_i $ the holderian constant and suppose that $ A_i \to 0 $, can one conclude with the blow-up technique that the sequence of the solutions $ u_i $ is compact ?

\section{Proofs of the results} 

\underbar {Proofs of the theorems:}

\smallskip

Without loss of generality, we can assume that $ \Omega= B_1(0) $ the unit ball centered on the origin.

\smallskip

We assume that:

$$ u_i \in  W_0^{1,1}(\Omega). $$

According to the work of Brezis-Merle, $ e^{ku_i} \in L^1 $ for all $ k >2 $ and the elliptic estimates and the Sobolev embedding gives:

$$ u_i \in W^{2, k}(\Omega) \cap C^{1, \epsilon}(\bar \Omega). $$

Here, $ G $ is the Green function of the Laplacian with Dirichlet condition on $ B_1(0) $. We have (in complex notation):
 
 $$ G(x, y)=\dfrac{1}{2 \pi} \log \dfrac{|1-\bar xy|}{|x-y|}, $$ 
 
we can write:

$$ v_i(x)=\int_{B_1(0)} G(x,y) V_i(y) e^{u_i(y)} dy , $$

Remark that, we can write:

$$ F(x)=\int_{B_1(0)}  -\dfrac{1}{2 \pi} \log |x-y| V_i(y) e^{u_i(y)} dy, $$

By a result in Gilbarg-Trudinger (cahpter 4, Newtonnian potential), we have $ F \in C^1(\bar \Omega) $

 Also, we have:

$$ G(x)=\int_{B_1(0)}  -\dfrac{1}{2 \pi} \log |1-\bar xy| V_i(y) e^{u_i(y)} dy, $$

For $ x $ near $ 0 $, $ G $ is smooth by the usuual differentability theorem. For $ x \not = 0 $, we can write:

$ K(x)= F(1/\bar x) $ which is $ C^1 $ by a result of Gilbarg-Trudinger. Combining the two last results, we have $ v_i $ is $ C^1 $ and by the maximum principle we have $ v_i \equiv u_i $. We use the fact that $ G $ is real to write $ \partial_x G = \bar  \partial_{x} \bar G = \bar  \partial_{\bar x}  G $ to have the dirivative of $ u_i $.

$$ G(x, y)=\dfrac{1}{2 \pi} \log \dfrac{|1-\bar xy|}{|x-y|}, $$

we can write:

$$ u_i(x)=\int_{B_1(0)} G(x,y) V_i(y) e^{u_i(y)} dy , $$

We can compute (in complex notation) $ \partial_x  G $ and $ \partial_x u_i $ :

 $$ \partial_x  G (x,y)=\dfrac{1-|y|^2}{(x-y)(x\bar y-1)}, $$
 
 $$ \partial_x u_i (x)=\int_{B_1(0)}  \partial_x  G (x,y) V_i(y) e^{u_i(y)} dy =\int_{B_1(0)} \dfrac{1-|y|^2}{(x-y)(x\bar y-1)} V_i(y) e^{u_i(y)} dy $$
  
we write,

$$ u_i(x_i)=\int_{\Omega} G(x_i,y) V_i(y) e^{u_i(y)} dx=\int_{\Omega-B(x_i,\delta_i/2)} G(x_i,y) V_i e^{u_i(y)} dy + \int_{B(x_i,\delta_i/2)} G(x_i,y) V_ie^{u_i(y)} dy  $$ 
 
According to the maximum principle,  the harmonic function $ G(x_i, .) $ on $ \Omega-B(x_i,\delta_i/2) $ take its maximum on the boundary of $ B(x_i,\delta_i/2) $, we can compute this maximum:

$$ G(x_i, y_i)=\dfrac{1}{2 \pi} \log \dfrac{|1-\bar x_iy_i|}{|x_i-y_i|} =\dfrac{1}{2\pi} \log \dfrac{|1-\bar x_i(x_i+\delta_i\theta_i)|}{|\delta_i/2| } =\dfrac{1}{2\pi} \log 2(|(1+|x_i|)+\bar x_i\theta_i|) < + \infty $$

with $ |\theta_i|= 1/2 $.

\smallskip

Thus,

$$ u_i(x_i) \leq C+  \int_{B(x_i,\delta_i/2)} G(x_i,y) V_ie^{u_i(y)} dy  \leq C + e^{u_i(x_i)}  \int_{B(x_i,\delta_i/2)} G(x_i,y) dy $$

Now, we compute $ \int_{B(x_i,\delta_i/2)} G(x_i,y) dy  $

we set in polar coordinates,

$$ y=x_i+\delta_i t \theta  $$

we find:

$$ \int_{B(x_i,\delta_i/2)} G(x_i,y) dy =\int_{B(x_i,\delta_i/2)} \dfrac{1}{2 \pi} \log \dfrac{|1-\bar x_iy|}{|x_i-y|} =\dfrac{1}{2 \pi}\int_0^{2\pi}\int_0^{1/2} \delta_i^2 \log \dfrac{|1-\bar x_i(x_i+\delta_i \theta)|}{(\delta_i t)/2} t dt d \theta = $$

$$ = \dfrac{1}{2 \pi}\int_0^{2\pi}\int_0^{1/2} \delta_i^2 (\log 2(|1+|x_i|+ t \bar x_i \theta|)-\log t) t dt d \theta \leq C \delta_i^2.$$

Thus,

$$ u_i(x_i) \leq C+ C \delta_i^2 e^{u_i(x_i)} , $$

which we can write, because $ u_i(x_i) \to +\infty $,

$$ u_i(x_i) \leq C' \delta_i^2 e^{u_i(x_i)}, $$

We can conclude that:

$$ u_i(x_i)+2 \log \delta_i \to + \infty. $$

Now, consider the following function :

$$ v_i(y)=u_i(x_i+\delta_i y)+2\log \delta_i, \quad y \in B(0,1/2) $$

The function satisfies all conditions of the Brezis-Merle hypothesis, we can conclude that, on each compact set:

$$ v_i \to -\infty $$

we can assume, without loss of generality that for $ 1/2 > \epsilon >0 $, we have:

$$ v_i \to -\infty, \quad y \in  B(0, 2\epsilon)-B(0,\epsilon), $$

\begin{Lemma} 

For all $ 1/4> \epsilon >0 $, we have:

$$ \sup_{B(x_i, (3/2)\delta_i\epsilon)-B(x_i, \delta_i \epsilon)} u_i \leq C_{\epsilon}. $$

\end{Lemma}

\underbar{Proof of the lemma}

\smallskip

Let $ t'_i $ and $ t_i $ the points of $ B(x_i, 2\delta_i\epsilon)-B(x_i, (1/2)\delta_i \epsilon) $ and $ B(x_i, (3/2)\delta_i\epsilon)-B(x_i, \delta_i \epsilon) $ respectively where $ u_i $ takes its maximum.

According to the Brezis-Merle work, we have:

$$ u_i(t'_i)+2\log \delta_i \to - \infty $$

We write,

$$ u_i(t_i)=\int_{\Omega} G(t_i,y) V_i(y) e^{u_i(y)} dx=\int_{\Omega-B(x_i,2\delta_i\epsilon)} G(t_i,y) V_i e^{u_i(y)} dy + $$

$$+ \int_{B(x_i,2\delta_i\epsilon)-B(x_i, (1/2)\delta_i\epsilon)} G(t_i,y) V_ie^{u_i(y)} dy+ $$

$$ +  \int_{B(x_i, (1/2)\delta_i\epsilon)} G(t_i,y) V_ie^{u_i(y)} dy $$ 

But, in the first and the third integrale, the point $ t_i $ is far from the singularity $ x_i $ and we know that the Green function is bounded. For the second integrale, after a change of variable, we can see that this integale is bounded by (we take the supremum in the annulus and use Brezis-Merle theorem)

$$ \delta_i^2 e^{u_i(t'_i)} \times  I_j $$

where $ I_j $ is a Jensen integrale (of the form $ \int_0^{1}  \int_0^{2\pi}  (\log (|1+|x_i|+ t \theta |)-\log |\theta_i-t\theta|) t dt d \theta $ which is bounded ). 

we conclude the lemma.

From the lemma, we see that far from the singularity the sequence is bounded, thus if we take the supremum on the set $ B_1(0)-B(x_i, \delta_i \epsilon) $ we can see that this supremum  is bounded and thus the sequence of functions is uniformly bounded or tends to infinity and we use the same arguments as for $ x_i $  to conclude that around this point and far from the singularity, the seqence is bounded.

The process will be finished , because, according to Brezis-Merle estimate, around each supremum constructed and tending to infinity, we have:

$$ \forall \,\, \epsilon >0, \,\,\, \limsup_{i \to + \infty} \int_{B(x_i, \delta_i\epsilon)} V_i e^{u_i}  dy \geq 4 \pi >0 .$$

Finaly, with this construction, we have a finite number of "exterior "blow-up points and outside the singularities the sequence is bounded uniformly, for example, in the case of one "exterior" blow-up point, we have:

$$ u_i(x_i) \to + \infty  $$

$$ \forall \,\, \epsilon >0, \,\, \sup_{B_1(0)-B(x_i, \delta_i \epsilon)} u_i \leq C_{\epsilon} $$

$$ \forall \,\, \epsilon >0, \,\,\, \limsup_{i \to + \infty} \int_{B(x_i, \delta_i\epsilon)} V_i e^{u_i}  dy \geq 4 \pi >0 .$$

$$ x_i \to x_0 \in \partial B_1(0). $$

We have the following lemma:

\begin{Lemma} Each $ \delta_i^k $ is of order $ d(x_i^k, \partial B_1(0)) $. Namely: there is a positive constant $ C >0 $ such that for $ \epsilon >0 $ small enough:

$$ \delta_i^k \leq d(x_i^k, \partial B_1(0)) \leq (2+\dfrac{C}{\epsilon})\delta_i^k. $$

\end{Lemma}

\underbar{Proof of the lemma}

Now, if we suppose that there is another "exterior" blow-up $ (t_i)_i $, we have, because $ (u_i)_i $ is uniformly bounded in a neighborhood of $ \partial B(x_i, \delta_i\epsilon) $, we have :

$$ d(t_i,\partial B(x_i, \delta_i\epsilon)) \geq \delta_i\epsilon $$

If we set,

$$ \delta'_i= d(t_i, \partial (B_1(0)-B(x_i, \delta_i\epsilon)))= \inf \{d(t_i,\partial B(x_i, \delta_i\epsilon)), d(t_i,\partial (B_1(0))) \} $$

then, $ \delta'_i  $ is of order $ d(t_i, \partial B_1(0)) $. To see this, we write:

$$ d(t_i,\partial B_1(0)) \leq d(t_i,\partial B(x_i, \delta_i\epsilon))+d( \partial B(x_i, \delta_i\epsilon), x_i) + d(x_i, \partial B_1(0)), $$

Thus,

$$ \dfrac{d(t_i,\partial B_1(0))}{d(t_i,\partial B(x_i, \delta_i\epsilon))} \leq 2+\dfrac{1}{\epsilon}, $$

Thus,

$$   \delta'_i \leq d(t_i,\partial B_1(0)) \leq \delta'_i (2+\dfrac{1}{\epsilon}). $$

Now, the general case follow by induction. We use the same argument for three, four,..., $ n $ blow-up points. 

We have, by induction and, here we use the fact that $ u_i $ is uniformly bounded outside a small ball centered at $ x_i^j, j=0,\ldots,k-1 $:

$$ \delta_i^j  \leq d(x_i^j, \partial B_1(0))\leq C_1 \delta_i^j, \,\, j=0,\ldots, k-1,  $$.

$$ d(x_i^k, \partial B(x_i^j, \delta_i^j \epsilon/2))\geq \epsilon \delta_i^j, \epsilon >0,\,\, j=0,\ldots, k-1,  $$.

and let's consider $ x_i^k $ such that:

$$ u_i(x_i^k)=\sup_{B_1(0)- \cup_{j=0}^{k-1}  B(x_i^j, \delta_i^j \epsilon)} u_i \to + \infty, $$

take,

$$ \delta_i^k= \inf \{d(x_i^k, \partial B_1(0)), d(x_i^k, \partial (B_1(0)-\cup_{j=0}^{k-1}  B(x_i^j, \delta_i^j \epsilon/2)) \}, $$

if, we have,

$$ \delta_i^k= d(x_i^k, \partial B(x_i^j, \delta_i^j \epsilon/2)), \,\, j\in\{0, \ldots,  k-1 \}.  $$ 

Then,

$$ \delta_i^k \leq d(x_i^k, \partial B_1(0)) \leq $$

$$ \leq d(x_i^k, \partial B(x_i^j, \delta_i^j \epsilon/2))+ d(\partial B(x_i^j, \delta_i^j \epsilon/2), x_i^j)+ d(x_i^j, \partial B_1(0))  $$

$$ \leq (2+\dfrac{C_1}{\epsilon})\delta_i^k. $$

To apply lemma 2.1 for $ m $ blow-up points, we use an induction:

We do directly the same approch for $ t_i $ as $ x_i $ by using directly the Green function of the unit ball.

\smallskip

If we look to the blow-up points, we can see, with this work that, after finite steps, the sequence will be bounded outside a finite number of balls , because of Brezis-Merle estimate (corollary 4 of Brezis-Merle's paper):

$$  \forall \,\, \epsilon >0, \,\,\, \limsup_{i \to + \infty} \int_{B(x_i^k, \delta_i^k\epsilon)} V_i e^{u_i}  dy \geq 4 \pi >0 .$$

Here, we can take the functions:

$$ u_i^k(y)=u_i(x_i^k+\delta_i^k y)+2\log \delta_i^k, $$

(By corollary 4 of Brezis-Merle's paper if $ \limsup_{i \to + \infty} \int_{B(x_i^k, \delta_i^k\epsilon)} V_i e^{u_i}  dy \leq 4\pi-\epsilon_0 < 4\pi $, then $ (u_i^k)^+ $ is locally uniformly bounded, which in contradiction with $ u_i^k(0) \to + \infty $). 

Finaly, we can say that, there is a finite number of sequences $ (x_i^k)_i, (\delta_i^k) , 0 \leq k \leq m $, such that:

$$  (x_i^0)_i \equiv (x_i)_i, \,\, \delta_i^0=\delta_i=d(x_i, \partial B_1(0)), $$

$$ (x_i^1)_i \equiv (t_i)_i, \,\, \delta_i^1=\delta'_i=d(t_i, \partial (B_1(0)-B(x_i, \delta_i \epsilon)), $$

and each $\delta_i^k $ is of order $ d(x_i^k, \partial B_1(0)) $.

and,

$$ u_i(x_i^k)=\sup_{B_1(0)- \cup_{j=0}^{k-1}  B(x_i^j, \delta_i^j \epsilon)} u_i \to + \infty, $$

$$ u_i(x_i^k) + 2\log \delta_i^k \to + \infty, $$

$$ \forall \,\, \epsilon >0, \,\, \sup_{B_1(0)- \cup_{j=0}^{m}  B(x_i^j, \delta_i^j \epsilon)} u_i \leq C_{\epsilon} $$

$$ \forall \,\, \epsilon >0, \,\,\, \limsup_{i \to + \infty} \int_{B(x_i^k, \delta_i^k\epsilon)} V_i e^{u_i}  dy \geq 4 \pi >0 .$$

\underbar {The work of YY.Li-I.Shafrir}

\smallskip

With the previous method, we have a finite number of "exterior" blow-up points (perhaps the same) and the sequences tend to the boundary. With the aid of proposition 1 of the paper of Li-Shafrir, we see that around each exterior blow-up, we have a finite number of "interior" blow-ups. Around, each exterior blow-up, we have after rescaling with $ \delta_i^k $, the same situation as around a fixed ball with positive radius.  If we assume:

$$ V_i \to V\,\, {\rm in } \,\, C^0(\bar B_1(0)), $$

then,

$$  \forall \,\, \epsilon >0, \,\,\, \limsup_{i \to + \infty} \int_{B(x_i^k, \delta_i^k\epsilon)} V_i e^{u_i}  dy = 8 \pi m_k, \,\, m_k \in  {\mathbb N}^*.$$

And, thus, we have the following convergence in the sense of distributions:

$$ \int_{B_1(0)} V_i e^{u_i}  dy  \to \int_{B_1(0)} V e^{u}  dy +\sum_{k=0}^m 8 \pi m_k' \delta_{x_0^k} ,\, m_k' \in  {\mathbb N}^*, \,\, x_0^k \in \partial B_1(0). $$

\underbar {Consequence 1: Proof of theorem 2}

\smallskip

Assume that:

$$ \int_{B_1(0)} V_i e^{u_i}  dy  \leq 4\pi, $$

Then, if the sequence blow-up, there is one and only one blow-up point and we have:

$$ u_i(x_i)=\sup_{B_1(0)} u_i \to + \infty, $$

$$ u_i(x_i) + 2\log \delta_i \to + \infty, $$

$$ \forall \,\, \epsilon >0, \,\, \sup_{B_1(0)- B(x_i, \delta_i \epsilon)} u_i \leq C_{\epsilon} $$

We set,

$$ r_i=e^{-u_i(x_i)/2}, $$

The blow-up function is locally bounded thus,

$$ r_i^2e^{u_i} \leq C\,\, {\rm on } \,\, B(x_i, 2r_i). $$

We write:

$$ u_i(x_i)=\int_{\Omega-B(x_i,\delta_i\epsilon)} G(x_i,y) V_i e^{u_i(y)} dy + \int_{B(x_i,\delta_i \epsilon )} G(x_i,y) V_ie^{u_i(y)} dy \leq  C_{\epsilon} + \int_{B(x_i,\delta_i \epsilon )} G(x_i,y) V_ie^{u_i(y)} dy $$ 
 
we have:

$$ \int_{B(x_i,\delta_i \epsilon )} G(x_i,y) V_ie^{u_i(y)} dy= \int_{B(x_i,\delta_i \epsilon )-B(x_i,2r_i)} G(x_i,y) V_ie^{u_i(y)} dy+ \int_{B(x_i, 2r_i )} G(x_i,y) V_ie^{u_i(y)} dy $$

We use the maximum principle on $ B(x_i,\delta_i \epsilon)-B(x_i,2r_i) $ and the explicit formula of $ G $ to prove that:

$$ G(x_i,y) \leq C+\dfrac{1}{2 \pi} \log \dfrac{\delta_i}{r_i} = C+ \dfrac{1}{4 \pi} (u_i(x_i)+2\log \delta_i). $$

On $ B(x_i, 2r_i) $ we use the fact that:

$$ r_i^2e^{u_i} \leq C $$

and the explicit formula for $ G $ to have:

$$ \int_{B(x_i, 2r_i )} G(x_i,y) V_ie^{u_i(y)} dy \leq C+\dfrac{1}{2 \pi} \log \dfrac{\delta_i}{r_i}\int_{B(x_i, 2r_i )}V_ie^{u_i(y)} dy. $$

We conclude that:

$$ u_i(x_i) \leq C+\dfrac{1}{2 \pi} \log \dfrac{\delta_i}{r_i}\int_{B(x_i, \delta_i \epsilon )}V_ie^{u_i(y)} dy. $$

which we can write as:

$$ u_i(x_i) \leq C+\dfrac{1}{4 \pi}(u_i(x_i)+2\log \delta_i)\int_{B(x_i, \delta_i \epsilon )}V_ie^{u_i(y)} dy. $$

Our hypothesis on the integrale of $ V_ie^{u_i} $ imply that:

$$ \log \delta_i  \geq -C, $$

in other words, we have uniformly,

$$  d(x_i, \partial B_1(0))= \delta_i  \geq e^{-C} >0. $$

this contredicts the fact that $ (x_i) $ tends to the boundary. The sequence $ (u_i) $ is bounded in this case.

We can see that the case:

$$ \int_{B_1(0)} V_i e^{u_i}  dy  \leq 4\pi, $$

is optimal, because Brezis-Merle have proved that, there is a counterexample of blow-up sequence with:

$$ \int_{B_1(0)} V_i e^{u_i}  dy  = 4\pi A >4\pi. $$

\underbar {Consequence 2: using a Pohozaev-type identity, proof of theorem 3}

\bigskip

By a conformal transformation, we can assume that our domain $ \Omega = B^+ $ is a half ball centered at the origin, $ B^+=\{x, |x|\leq 1,  x_1 \geq 0 \} $. In this case the normal at the boundary is $ \nu=(-1,0) $ and $ u_i(0, x_2) \equiv 0 $. Also, we set $ x_i $ the blow-up point and $ x_i^2=(0,x_i^2) $ and $ x_i^1=(x_i^1,0) $ respectevely the second and the first part of $ x_i $. Let $ \partial B^+ $ the part of the boundary for which $ u_i $ and its derivatives are uniformly bounded and thus converge to the corresponding function.

\smallskip

\underbar {The case of one blow-up point:}

\smallskip

\begin{Theorem} If  $ V_i $ is s-Holderian with $ 1/2 < s  \leq 1$ and,

$$ \int_{\Omega}V_i e^{u_i} dy \leq 16\pi-\epsilon , \,\, \epsilon >0, $$

we have :
 
 $$ V_i(x_i)\int_{\Omega} e^{u_i} dy- V(0)\int_{\Omega} e^{u} dy  =o(1)$$
 
 which means that there is no blow-up points.

\end{Theorem}

\smallskip

\underbar{Proof of the theorem}

\smallskip

The Pohozaev identity gives us the following formula:

$$ \int_{\Omega} <(x-x_2^i) |\nabla u_i>(-\Delta u_i) dy =\int_{\Omega} <(x-x_2^i) |\nabla u_i> V_i e^{u_i} dy = A_i $$

$$ A_i= \int_{\partial B^+} <(x-x_2^i) |\nabla u_i> <\nu |\nabla u_i> d\sigma + \int_{\partial B^+} <(x-x_2^i) |\nu>|\nabla u_i|^2d\sigma $$

We can write it as:

$$ \int_{\Omega} <(x-x_2^i) |\nabla u_i> (V_i-V_i(x_i)) e^{u_i} dy =A_i+  V_i(x_i) \int_{\Omega} <(x-x_2^i) |\nabla u_i> e^{u_i} dy = $$

$$ = A_i+ V_i(x_i)\int_{\Omega} <(x-x_2^i) |\nabla (e^{u_i})> dy  $$

And, if we integrate by part the second term, we have (because $ x_1=0 $ on the boundary and $ \nu_2=0 $):

$$ \int_{\Omega} <(x-x_2^i) |\nabla u_i> (V_i-V_i(x_i)) e^{u_i} dy =-2V_i(x_i)\int_{\Omega} e^{u_i} dy +B_i $$

where $ B_i $ is,

$$ B_i=V_i(x_i)\int_{\partial B^+}  <(x-x_2^i) |\nu> e^{u_i} dy $$

applying the same procedure to $ u $, we can write:

$$ -2V_i(x_i)\int_{\Omega} e^{u_i} dy+2V(0)\int_{\Omega} e^{u} dy = \int_{\Omega} <(x-x_2^i) |\nabla u_i> (V_i-V_i(x_i)) e^{u_i} dy-\int_{\Omega} <(x-x_2^i) |\nabla u> (V-V(0)) e^{u} dy+ $$

$$ + (A_i-A)+ (B_i-B), $$

where $ A $ and $ B $ are,

$$ A= \int_{\partial B^+} <(x-x_2^i) |\nabla u> <\nu |\nabla u> d\sigma + \int_{\partial B^+} <(x-x_2^i) |\nu>|\nabla u|^2d\sigma $$

$$ B=V(0)\int_{\partial B^+}  <(x-x_2^i) |\nu> e^{u} dy $$

and, because of the uniform convergence of $ u_i $ and its derivatives on $ \partial B^+ $, we have:

$$ A_i-A = o(1) \,\,\, {\rm and}\,\,\, B_i-B = o(1) $$

which we can write as:

$$ V_i(x_i)\int_{\Omega} e^{u_i} dy- V(0)\int_{\Omega} e^{u} dy = \int_{\Omega} <(x-x_2^i) |\nabla (u_i-u)> (V_i-V_i(x_i)) e^{u_i} dy +  $$

$$ +  \int_{\Omega} <(x-x_2^i) |\nabla u> (V_i-V_i(x_i)) (e^{u_i}-e^u) dy + $$

$$  + \int_{\Omega} <(x-x_2^i) |\nabla u> (V_i-V_i(x_i)- (V-V(0))) e^{u} dy + o(1) $$ 

We can write the second term as:

 $$ \int_{\Omega} <(x-x_2^i) |\nabla u> (V_i-V_i(x_i)) (e^{u_i}-e^u) dy = \int_{\Omega-B(0,\epsilon) } <(x-x_2^i) |\nabla u> (V_i-V_i(x_i)) (e^{u_i}-e^u) dy + $$
 
 $$ + \int_{B(0,\epsilon)} <(x-x_2^i) |\nabla u> (V_i-V_i(x_i)) (e^{u_i}-e^u) dy = o(1), $$

because of the uniform convergence of $ u_i $ to $ u $ outside a region which contain the blow-up and the uniform convergence of $ V_i $. For the third integral we have the same result:

$$ \int_{\Omega} <(x-x_2^i) |\nabla u> (V_i-V_i(x_i)- (V-V(0))) e^{u} dy =o(1), $$

because of the uniform convergence of $ V_i  $ to $ V $.

Now, we look to the first integral:

$$ \int_{\Omega} <(x-x_2^i) |\nabla (u_i-u)> (V_i-V_i(x_i)) e^{u_i} dy,  $$

we can write it as:

$$ \int_{\Omega} <(x-x_2^i) |\nabla (u_i-u)> (V_i-V_i(x_i)) e^{u_i} dy =\int_{\Omega} <(x-x_i) |\nabla (u_i-u)> (V_i-V_i(x_i)) e^{u_i} dy + $$

$$ + \int_{\Omega} <x_1^i |\nabla (u_i-u)> (V_i-V_i(x_i)) e^{u_i} dy, $$

Thus, we have proved by using the Pohozaev identity  the following equality, in the case of one blow-up (the term $ o(1) $):

$$ \int_{\Omega} <(x-x_i) |\nabla (u_i-u)> (V_i-V_i(x_i)) e^{u_i} dy + $$

$$ + \int_{\Omega} <x_1^i |\nabla (u_i-u)> (V_i-V_i(x_i)) e^{u_i} dy = $$

$$ = 2V_i(x_i)\int_{\Omega} e^{u_i} dy- 2V(0)\int_{\Omega} e^{u} dy  +o(1)$$

We can see, because of the uniform boundedness of $ u_i $ outside $ B(x_i, \delta_i \epsilon) $ and the fact that :

$$ ||\nabla (u_i-u) ||_1 = o(1), $$

it is sufficient to look to the integral on $ B(x_i, \delta_i \epsilon) $.

\smallskip

Assume that we are in the case of one blow-up, it must be $ (x_i) $ and isolated, we can write the following inequality as a consequence of YY.Li-I.Shafrir result:

$$ u_i(x)+ 2\log |x-x_i| \leq C,  $$

We use this fact and the fact that $ V_i $ is s-holderian to have that, on $ B(x_i, \delta_i \epsilon) $,

$$ |(x-x_i)(V_i-V_i(x_i)) e^{u_i}| \leq \dfrac{C}{ |x-x_i|^{1-s}} \in L^{(2-\epsilon')/(1-s)}, \,\, \forall \,\, \epsilon' >0, $$

and, we use the fact that:

$$ ||\nabla (u_i-u) ||_q = o(1), \,\, \forall \,\, 1 \leq q < 2$$

to conclude by the Holder inequality that for $ 0 < s \leq 1 $ :

$$ \int_{B(x_i, \delta_i \epsilon)} <(x-x_i) |\nabla (u_i-u)> (V_i-V_i(x_i)) e^{u_i} dy = o(1), $$

For the other integral, namely:

$$  \int_{B(x_i, \delta_i \epsilon)} <x_1^i |\nabla (u_i-u)> (V_i-V_i(x_i)) e^{u_i} dy , $$

We use the fact that, because our domain is a half ball, and the $ \sup+\inf  $ inequality to have:

$$ x_1^i= \delta_i, $$

$$ u_i(x)+4 \log \delta_i \leq C $$

and,

$$ e^{ (s/2) u_i(x)} \leq |x-x_i|^{-s},  $$

$$ |V_i-V_i(x_i)|\leq |x-x_i|^s, $$

Finaly, we have:

$$  | \int_{B(x_i, \delta_i \epsilon)} <x_1^i |\nabla (u_i-u)> (V_i-V_i(x_i)) e^{u_i} dy |\leq C \int_{B(x_i, \delta_i \epsilon)}  |\nabla (u_i-u)| e^{((3/4)-(s/2))u_i}, $$

But in the second member, for $ 1/2 < s \leq 1 $, we have $ q_s=1/(3/4-s/2) >2 $ and thus $ q_s' <2 $ and,

$$ e^{((3/4)-(s/2))u_i} \in L^{q_s} $$

$$  ||\nabla (u_i-u) ||_{q_s'}  = o(1), \,\, \forall \,\, 1 \leq q_s' < 2, $$

one conclude that:

$$  \int_{B(x_i, \delta_i \epsilon)} <x_1^i |\nabla (u_i-u)> (V_i-V_i(x_i)) e^{u_i} dy =o(1) $$

 Finaly, with this method, we conclude that, in the case of one blow-up point and $ V_i $ is s-Holderian with $ 1/2 < s  \leq 1$ :

 $$ V_i(x_i)\int_{\Omega} e^{u_i} dy- V(0)\int_{\Omega} e^{u} dy  =o(1)$$
 
 which means that there is no blow-up, which is a contradiction.

 Finaly, for one blow-up point and $ V_i  $ is s-Holderian with $ 1/2 < s  \leq 1$, the sequence $ (u_i) $ is uniformly bounded on $ \Omega $. 
 
\bigskip

\underbar {The case of two blow-up points:}

\smallskip

\begin{Theorem} If  $ V_i $ is s-Holderian with $ 1/2 < s  \leq 1$ and,

$$ \int_{\Omega}V_i e^{u_i} dy \leq 24\pi-\epsilon , \,\, \epsilon >0, $$

we have :
 
 $$ V_i(x_i)\int_{\Omega} e^{u_i} dy- V(0)\int_{\Omega} e^{u} dy  =o(1)$$
 
 which means that there is no blow-up points.

\end{Theorem}

\smallskip

\underbar{Proof of the Theorem}

\bigskip
 
\underbar {The case of two "interior" blow-up points:}

\smallskip
 
As in the previous case, we assume that  $  \Omega = B^+ $ is the half ball. We have two "interior" blow-up points $ x_i $ and $ y_i $:

$$ |y_i-x_i|\leq \delta_i \epsilon, $$

We use a Pohozaev type identity:

$$ \int_{\Omega} <(x-x_2^i) |\nabla u_i>(-\Delta u_i) dy =\int_{\Omega} <(x-x_2^i) |\nabla u_i> V_i e^{u_i} dy = A_i $$

with $ A_i $ the regular part of the identity (on which the uniform convergence holds).

$$ A_i= \int_{\partial B^+} <(x-x_2^i) |\nabla u_i> <\nu |\nabla u_i> d\sigma + \int_{\partial B^+} <(x-x_2^i) |\nu>|\nabla u_i|^2d\sigma $$

We divide our domain in two domain $ \Omega_1^i $ and $ \Omega_2^i $ such that:

$$ \Omega_1^i= \{x, |x-x_i|\leq |x-y_i| \}, \,\, \Omega_2^i= \{x, |x-x_i|\geq |x-y_i| \}. $$

We set,

$$ D_i= \{x, |x-x_i|=|x-y_i| \}. $$

We write:

$$ A_i=\int_{\Omega_1^i} <(x-x_2^i) |\nabla u_i> (V_i -V_i(x_i))e^{u_i} dy +\int_{\Omega_2^i} <(x-x_2^i) |\nabla u_i> (V_i-V_i(y_i)) e^{u_i} dy + $$

$$ + V_i(x_i)\int_{\Omega_1^i} <(x-x_2^i) |\nabla u_i>  e^{u_i} dy+ V_i(y_i)\int_{\Omega_2^i} <(x-x_2^i) |\nabla u_i>  e^{u_i} dy. $$

As for the case of one blow-up point, it is sufficient to consider terms  which contain  the difference $\nabla (u_i-u) $.

We can write the last addition as (after using $ \nabla (u_i-u) $) :

$$ \left ( V_i(x_i) \int_{\Omega} <(x-x_2^i) |\nabla u_i> e^{u_i} dy - \int_{\Omega} <(x-x_2^i) |\nabla u> e^{u} dy \right ) + $$

$$ + ( V_i(y_i)-V_i(x_i))\int_{\Omega_2^i} <(x-x_2^i) |\nabla (u_i-u)> e^{u_i} dy. $$

First of all, we consider the term (which equal, after integration by part to ):

$$ V_i(x_i) \int_{\Omega} <(x-x_2^i) |\nabla u_i> e^{u_i} dy- \int_{\Omega} <(x-x_2^i) |\nabla u> e^{u} dy = $$

$$ =-2V_i(x_i) \int_{\Omega} e^{u_i} dy +2V(0) \int_{\Omega} e^{u} dy + (B_i-B) $$

with the same notation for $ B_i $ and $ B $ as for the previous case.

\bigskip

\underbar {Case 1:} suppose that, $ |x-y_i|\geq |x_i-y_i| $, 

thus

$$ |V_i(x_i)-V_i(y_i)|\leq |x_i-y_i|^s \leq |x-y_i|^s $$

Thus,

$$|( V_i(y_i)-V_i(x_i))\int_{\Omega_2^i \cap \{x, |x-x_i|\geq |x-y_i| \}.} <(x-x_2^i) |\nabla (u_i-u)> e^{u_i} dy | \leq \int_{\Omega_2^i} |x-y_i|^{1+s} |\nabla (u_i-u)| e^{u_i} dy+ $$

$$ + |y_2^i-x_2^i|\int_{\Omega_2^i} |x-y_i|^{s} |\nabla (u_i-u)| e^{u_i} dy+|y_1^i|\int_{\Omega_2^i} |x-y_i|^{s} |\nabla (u_i-u)| e^{u_i} dy $$

But,

$$ |y_i-x_i|\leq \delta_i \epsilon, $$

$$ x_1^i = \delta_i $$

we use the same method (with the $ \sup+\inf $ inequality) to prove that for $  1 \geq s >1/2 $ the two integrals converges to $ 0 $.

\underbar {Case 2:} suppose that, $ |x-y_i|\leq |x_i-y_i| $, 

We do integration by parts, we have one part on $ D_i $ and the other one on  the circle with center $ y_i $. In fact, we have intersection of convex 2-dimensional domains, which is convex (Lipschitz domain or a domain with finite nimber of singularity) and thus one can apply the Stokes or Green-Riemann formulas.

$$ (V_i(y_i)-V_i(x_i))\int_{\Omega_2^i \cap \{x, |x-y_i|\leq |x_i-y_i| \}.} <(x-x_2^i) |\nabla (e^{u_i})> dy = $$

$$ =( V_i(y_i)-V_i(x_i))\int_{D_i \cap \{x, |x-y_i|\leq |x_i-y_i| \}.} <(x-x_2^i) |\nu> e^{u_i} dy + $$

$$ + ( V_i(y_i)-V_i(x_i))\int_{ \{x, |x-y_i|= |x_i-y_i| \}  \cap \{x, |x-y_i|\leq |x-x_i| \}} <(x-x_2^i) |\nu> e^{u_i} dy +  $$

$$ + 2( V_i(y_i)-V_i(x_i))\int_{ \{x, |x-y_i|\leq |x_i-y_i| \}} e^{u_i} dy $$

We set:

$$ I_1=( V_i(y_i)-V_i(x_i))\int_{D_i \cap \{x, |x-y_i|\leq |x_i-y_i| \}.} <(x-x_2^i) |\nu> e^{u_i} dy, $$

$$ I_2= ( V_i(y_i)-V_i(x_i))\int_{ \{x, |x-x_i|= |x_i-y_i| \}  \cap \{x, |x-y_i|\leq |x-x_i| \}} <(x-x_2^i) |\nu> e^{u_i} dy $$

\begin{Lemma} We have:

$$ I_1=o(1), $$

and,

$$ I_2=o(1). $$

\end{Lemma}

\underbar {Proof of the lemma} 

\bigskip

For $ I_1 $, we have:

$$ |V_i(x_i)-V_i(y_i)|\leq 2C |x-y_i|^s ,$$

$$ |I_1| \leq C\int_{D_i \cap \{x, |x-y_i|\leq |x_i-y_i| \}.} |<(x-y^i) |\nu> ||x-y_i|^s e^{u_i}+ $$

$$+ |x_2^i-y_2^i| \int_{D_i \cap \{x, |x-y_i|\leq |x_i-y_i| \}.} |x-y_i|^s e^{u_i} dy + $$

$$ + |y_1^i| \int_{D_i \cap \{x, |x-y_i|\leq |x_i-y_i| \}.} |x-y_i|^s e^{u_i} dy $$

But, 

$$ x_1^i= \delta_i, $$

$$ |y_i-x_i|\leq \delta_i \epsilon, $$

$$ u_i(x)+4 \log \delta_i \leq C $$

and,

$$ e^{ (3/4) u_i(x)} \leq |x-y_i|^{-3/2},  $$

Thus,

$$ |I_1| \leq\int_{D_i \cap \{x, |x-y_i|\leq |x_i-y_i| \}.}  | x-y_i|^{s-1}+ $$

$$+ C\int_{D_i \cap \{x, |x-y_i|\leq |x_i-y_i| \}.} |x-y_i|^{(-3/2)+s} dy, $$

If we set $ t_0=(x_i+y_i)/2 $, we have on one part of $ D_i $:

$$ |x-t_0|\leq |x-y_i| = |x-x_i|\leq |x_i-y_i|, $$

by a change of variable $ u=x-t_0 $ on the line $ D_i $, we can compute the two last integrals directly, to have, for $ 1 \geq s >1/2 $:

$$  |I_1| \leq C(| x_i-y_i|^{s} + |x_i-y_i|^{s-(1/2)}) = o(1), $$

For $ I_2 $ we have:

$$ I_2= ( V_i(y_i)-V_i(x_i))\int_{ \{x, |x-y_i|= |x_i-y_i| \}  \cap \{x, |x-x_i|\leq |x_i-y_i| \}} <(x-x_2^i) |\nu> e^{u_i} dy $$

and,

$$ |V_i(x_i)-V_i(y_i)|\leq 2C |x-y_i|^s ,$$

$$ |I_2| \leq C\int_{\{x, |x-y_i|= |x_i-y_i| \} \cap \{x, |x-y_i|\leq |x-x_i| \}.} |<(x-y_i) |\nu> ||x-y_i|^s e^{u_i}+ $$

$$+ |x_2^i-y_2^i| \int_{\{x, |x-y_i|= |x_i-y_i| \}  \cap \{x, |x-y_i|\leq |x-x_i| \}.} |x-y_i|^s e^{u_i} dy + $$

$$ + |y_1^i| \int_{\{x, |x-y_i|= |x_i-y_i| \}  \cap \{x, |x-x_i|\leq |x-x_i| \}.} |x-y_i|^s e^{u_i} dy $$

with the same method as for $ I_1 $ we have:

$$ |I_2| \leq C\int_{\{x, |x-y_i|= |x_i-y_i| \} \cap \{x, |x-y_i|\leq |x-x_i| \}.}  | x-y_i|^{s-1}+ $$

$$+\int_{\{x, |x-y_i|= |x_i-y_i| \}  \cap \{x, |x-y_i|\leq |x-x_i| \}.} |x-y_i|^{-(3/2)+s} dy, $$

Finaly, we have:

$$  |I_2| \leq C(| x_i-y_i|^{s} + |x_i-y_i|^{s-(1/2)}) = o(1), $$
 
\underbar {The case of two "exterior" blow-up points:}

\bigskip

Let $ (x_i)_i $ and $ (t_i)_i $ two sequences of "exterior" blow-up points. If $ d(x_i,t_i)=O( \delta_i) $ or $  d(x_i,t_i)=O( \delta_i') $ then we use  the same technique as for two interior blow-up with the Pohozaev identity. In this case the $  \sup+  \inf  $ inequality holds, because $ d(x_i, t_i) $ is of order $  \delta_i $ or $  \delta_i' $. Assume that:

$$  \dfrac{d(x_i,t_i)}{\delta_i}\to + \infty \,\,\, {\rm and } \,\,\, \dfrac{d(x_i,t_i)}{\delta_i'}\to + \infty $$

In this case, we assume that, we are on the half ball.(In fact one consider the intersection of disks  and half plane which is convex, and we take its image by the conformal map. In this case we have a domain on which we can apply the Stokes formula). By a conformal transformation, $ f $, we can assume that our two sequences  are on the unit ball. First of all, we use the Pohozaev identity on the half ball as for the previous cases, but our domain change, we have one part is vertical, the second part is a part of the boundary of the unit ball, in which the sequences $ (u_i) $ and $ (\partial u_i)_i $ are uniformly bounded and converge to the corresponding function, and the third part of boundary, is a regular curve $ D_i' $ such that its image by $ f $ is the mediatrice $ D_i $ of  the segment $ (x_i, t_i) $. In the Pohozaev identity, we have a terms of type:

$$ \int_{D_i'} <(x-x_2^i) |\nabla u_i> <\nu |\nabla u_i> d\sigma + \int_{D_i'} <(x-x_2^i) |\nu>|\nabla u_i|^2d\sigma $$

But if we integrate on the rest of the domain and if we use the Pohozeav identity on this second domain and we replace $ x_2^i $ by $ t_2^i $, the integral on $ D_i' $ is :

$$ - \int_{D_i'} <(x-t_2^i) |\nabla u_i> <\nu |\nabla u_i> d\sigma - \int_{D_i'} <(x-t_2^i) |\nu>|\nabla u_i|^2d\sigma $$

If, we add the two integral, we find:

$$\int_{D_i'} <(x_2^i-t_2^i) |\nabla u_i> <\nu |\nabla u_i> d\sigma + \int_{D_i'} <(x_2^i-t_2^i) |\nu>|\nabla u_i|^2d\sigma $$

We have the same techniques as for the previous cases ("interior" blow-up), except the fact that here, we use the Pohozaev identity on two differents domains which the union is our half ball.( In fact the image by a conformal map of a convex domain, intersection of two disks and a half plane. (Intersection of the unit disk and $ D(x_0,\epsilon) $ and the mediatrice of $[x_i, t_i] $, here, $ x_0 $ is the "blow-up" point, and we apply the conformal map $ f $)). And apply the Green-Riemann theorem for smooth domains with finite number of singular points on the boundary. Or directly the Stokes theorem with the fact that we have a Lipschitz domain because it is the image of a Lipschitz domain by a conformal map (Hofmann-Mitrea-Taylor).

Remark that, here, because we have the two conditions on $ d(x_i,t_i) $ and $ \delta_i, \delta_i' $, the mediatrice of $ [x_i, t_i] $ is close to a fixed segment $ [0, x_0] $. (use angles for this fact).

\bigskip

To conclude, we must show that this last integral is close to $ 0 $ as $ i $ tends to $ + \infty $. By a conformal map, it is sufficient to prove that the corresponding integral on the unit ball on $ D_i $ tends to $ 0 $. Without loss of generality, we can assume here that we work on the unit ball (for this integral).

\bigskip

On the unit ball, with the Dirichlet condition, the Green function is (in complex notation) :

$$ G(x, y)=\dfrac{1}{2 \pi} \log \dfrac{|1-\bar xy|}{|x-y|}, $$

we can write:

$$ u_i(x)=\int_{B_1(0)} G(x,y) V_i(y) e^{u_i(y)} dy , $$

We can compute (in complex notation) $ \partial_x  G $ and $ \partial_x u_i $ :

 $$ \partial_x  G (x,y)=\dfrac{1-|y|^2}{(x-y)(x\bar y-1)}, $$
 
 $$ \partial_x u_i (x)=\int_{B_1(0)}  \partial_x  G (x,y) V_i(y) e^{u_i(y)} dy =\int_{B_1(0)} \dfrac{1-|y|^2}{(x-y)(x\bar y-1)} V_i(y) e^{u_i(y)} dy $$

Let $ t_0^i= (x_i+t_i)/2 $. We assume that $ |x-t_0^i|\leq 1-\epsilon $ and   $ |t_0^i|\geq 1-(\epsilon/2) $.

\begin{Proposition} 1) For  $ ((1/2)+\tilde \epsilon)|x_i-t_i| \leq |x-t_0^i|\leq 1-\epsilon $ we have,

$$  |\partial_x u_i (x)|\leq C'+ C\dfrac{\delta_i}{|x_i-t_i|}\dfrac{1}{|x-t_0^i|}=C'+\dfrac{o(1)}{|x-t_0^i|}. $$

2) For  $ |x-t_0^i| \leq((1/2)-\tilde\epsilon)|x_i-t_i| $ we have,

$$  |\partial_x u_i (x)|\leq  C'+C\dfrac{\delta_i}{|x_i-t_i|}\dfrac{1}{|x_i-t_0^i|}=C'+\dfrac{o(1)}{|x_i-t_0^i|}. $$

with $ o(1) \to 0  $ as $ i \to + \infty $ .

\bigskip

3) For  $ ((1/2)-\tilde\epsilon)|x_i-t_i|\leq |x-t_0^i| \leq((1/2)+\tilde\epsilon)|x_i-t_i| $ we have,

$$ |x_i-t_i||\nabla u_i|_{L^{\infty}(D_i\cap \{ ((1/2)-\tilde\epsilon)|x_i-t_i|\leq |x-t_0^i| \leq((1/2)+\tilde\epsilon)|x_i-t_i| \}} \leq C. $$

\end{Proposition}

\underbar {Proof of the proposition:}

\bigskip

To estimate $ \partial_x u_i $ on $ D_i $, we divide the last integral in three parts:

$$ \partial_x u_i (x)=\int_{B_1(0)-(B(x_i,\delta_i\epsilon)\cup B(t_i, \delta_i'\epsilon))} \dfrac{1-|y|^2}{(x-y)(x\bar y-1)} V_i(y) e^{u_i(y)} dy + $$ 

$$ + \int_{B(x_i,\delta_i\epsilon)} \dfrac{1-|y|^2}{(x-y)(x\bar y-1)} V_i(y) e^{u_i(y)} dy + $$

$$ + \int_{B(t_i,\delta_i'\epsilon)} \dfrac{1-|y|^2}{(x-y)(x\bar y-1)} V_i(y) e^{u_i(y)} dy $$

Let us set:

$$ I_1=\int_{B_1(0)-(B(x_i,\delta_i\epsilon)\cup B(t_i, \delta_i'\epsilon))} \dfrac{1-|y|^2}{(x-y)(x\bar y-1)} V_i(y) e^{u_i(y)} dy $$

$$ I_2=\int_{B(x_i,\delta_i\epsilon)} \dfrac{1-|y|^2}{(x-y)(x\bar y-1)} V_i(y) e^{u_i(y)} dy,
$$

$$ I_3= \int_{B(t_i,\delta_i'\epsilon)} \dfrac{1-|y|^2}{(x-y)(x\bar y-1)} V_i(y) e^{u_i(y)} dy $$

For the first integral, because $ u_i \leq C $ on $ B_1(0)-(B(x_i,\delta_i\epsilon)\cup B(t_i, \delta_i'\epsilon)) $, we have:

$$ |I_1|\leq  C\int_{B_1(0)} \dfrac{1-|y|^2}{|x-y||x\bar y-1|} dy, $$ 

But, $ 1 \geq |x| = |x-t_0^i+t_0^i|\geq |t_0^i|-|x-t_0^i|\geq 1-(\epsilon/2)-(1-\epsilon)=\epsilon/2 $, thus, we can write:

$$ |I_1|\leq  C\int_{B_1(0)} \dfrac{1-|y|^2}{|x-y||x||\bar y-1/x|} dy, $$

and, we use the fact that:

$$ |\bar y-1/x| \geq ||\bar y|-1/|x||\geq |1/|x|-|y||\geq (1-|y|), $$

To have:

$$ |\partial_x u_i (x)|\leq |I_2|+|I_3|+  C\int_{B_1(0)} \dfrac{1+|y|}{|x-y|} dy =|I_2|+|I_3|+ C', $$

Now, we look to the second and third integrals, it is sufficient to consider the first one :

$$I_2= \int_{B(x_i,\delta_i\epsilon)} \dfrac{1-|y|^2}{(x-y)(x\bar y-1)} V_i(y) e^{u_i(y)} dy $$

\underbar {Case 1: $ ((1/2)+\tilde \epsilon)|x_i-t_i| \leq |x-t_0^i|\leq 1-\epsilon $:}

\bigskip

In this case we have:

$$ 1-|y|^2=1-|x_i+\delta_i z|^2=\delta_i(2+o(1)), $$

and,

$$ |x-y|=|x-t_0^i+t_0^i-y_i-\delta_i z|\geq (\tilde \epsilon/2)|x_i-t_i|, $$

and,

$$ |x\bar y-1|=|((x-t_0^i+t_0^i-x_i)+x_i)(\bar x_i+\delta_i \bar z)-1|\geq (\tilde \epsilon/2)|x-t_0^i|,  $$

Thus,

$$  |\partial_x u_i (x)|\leq C'+ C\dfrac{\delta_i}{|x_i-t_i|}\dfrac{1}{|x-t_0^i|}=C'+\dfrac{o(1)}{|x-t_0^i|}. $$

with, $ o(1) \to 0  $ as $ i \to + \infty $ .

\bigskip 

\underbar {Case 2: $ |x-t_0^i| \leq((1/2)-\tilde\epsilon)|x_i-t_i| $:}

\bigskip

In this case, we have:

$$ 1-|y|^2=1-|x_i+\delta_i z|^2=\delta_i(2+o(1)), $$

and,

$$ |x-y|=|x-t_0^i+t_0^i-y_i-\delta_i z|\geq (\tilde \epsilon/2)|x_i-t_i|, $$

and,

$$ |x\bar y-1|=|((x-t_0^i+t_0^i-x_i)+x_i)(\bar x_i+\delta_i \bar z)-1|\geq (\tilde \epsilon/2)|x_i-t_0^i|,  $$

Thus,

$$  |\partial_x u_i (x)|\leq  C'+C\dfrac{\delta_i}{|x_i-t_i|}\dfrac{1}{|x_i-t_0^i|}=C'+\dfrac{o(1)}{|x_i-t_0^i|}. $$

with, $ o(1) \to 0  $ as $ i \to + \infty $ .

\bigskip 

\underbar {Case 3: $ ((1/2)-\tilde\epsilon)|x_i-t_i|\leq |x-t_0^i| \leq((1/2)+\tilde\epsilon)|x_i-t_i| $:}

\bigskip

Let $ \tilde t_0^i $ the point of $ D_i $ such that $ |\tilde t_0^i-t_0^i|=1/2(|x_i-t_i|) $. We use the fact that the function:

$$ v_i(t)=u_i(\tilde t_0^i + (|x_i-t_0^i|/4) t), $$

is uniformly bounded for $ |t| \leq 1 $ and is a solution of PDE which is uniformly bounded on $ |t| \leq 1 $. By the elliptic estimates we have:

$$ |x_i-t_i||\nabla u_i|_{L^{\infty}(D_i\cap \{ ((1/2)-\tilde\epsilon)|x_i-t_i|\leq |x-t_0^i| \leq((1/2)+\tilde\epsilon)|x_i-t_i| \}} \leq C. $$

Thus, we use the previous cases to compute the following integral:

$$\int_{D_i} <(x_2^i-t_2^i) |\nabla u_i> <\nu |\nabla u_i> d\sigma + \int_{D_i} <(x_2^i-t_2^i) |\nu>|\nabla u_i|^2d\sigma = o(1) $$

and, thus,

$$\int_{D_i'} <(x_2^i-t_2^i) |\nabla u_i> <\nu |\nabla u_i> d\sigma + \int_{D_i'} <(x_2^i-t_2^i) |\nu>|\nabla u_i|^2d\sigma = o(1) $$

here, we used the previous estimates with $ i\to +\infty $ and $ \tilde \epsilon \to 0 $ (for the previous case 3).


\bigskip 

\begin{thebibliography}{99} 

\bibitem{1}{T. Aubin. Some Nonlinear Problems in Riemannian Geometry. Springer-Verlag 1998 }


\bibitem{2}{ S.S Bahoura. Majorations du type $ \sup u \times \inf u \leq c $ pour l'\'equation de la courbure scalaire sur un ouvert de $ {\mathbb R}^n, n\geq 3 $. J. Math. Pures. Appl.(9) 83 2004 no, 9, 1109-1150.}

\bibitem{3}{S.S. Bahoura. Harnack inequalities for Yamabe type equations.  Bull. Sci. Math.  133  (2009),  no. 8, 875-892}

\bibitem{4}{S.S. Bahoura. Lower bounds for sup+inf and sup $ \times $ inf and an extension of Chen-Lin result in dimension 3.  Acta Math. Sci. Ser. B Engl. Ed.  28  (2008),  no. 4, 749-758}

\bibitem{5}{S.S. Bahoura. Estimations uniformes pour l'�quation de Yamabe en dimensions 5 et 6. J. Funct. Anal.  242  (2007),  no. 2, 550-562.}

\bibitem{6}{C. Bandle. Isoperimetric inequalities and Applications. Pitman. 1980.}

\bibitem{6}{H. Brezis, YY. Li , I. Shafrir. A sup+inf inequality for some
nonlinear elliptic equations involving exponential
nonlinearities. J.Funct.Anal.115 (1993) 344-358.
}

\bibitem{7}{H.Brezis and F.Merle, Uniform estimates and blow-up bihavior for solutions of $ -\Delta u=Ve^u $ in two dimensions, Commun Partial Differential Equations 16 (1991), 1223-1253.
}

\bibitem{8}{L. Caffarelli, B. Gidas, J. Spruck. Asymptotic symmetry and local
behavior of semilinear elliptic equations with critical Sobolev
growth. Comm. Pure Appl. Math. 37 (1984) 369-402.
}

\bibitem{9}{W. Chen, C. Li. A priori Estimates for solutions to Nonlinear Elliptic Equations. Arch. Rational. Mech. Anal. 122 (1993) 145-157.}


\bibitem{10}{C-C.Chen, C-S. Lin. Estimates of the conformal scalar curvature
equation via the method of moving planes. Comm. Pure
Appl. Math. L(1997) 0971-1017.}

\bibitem{11}{C-C.Chen, C-S. Lin. A sharp sup+inf inequality for a nonlinear elliptic equation in ${\mathbb R}^2$.
Commun. Anal. Geom. 6, No.1, 1-19 (1998).}

\bibitem{12}{C-C.Chen, C-S. Lin. Sharp estimates for solutions of multi-bubbles in compact Riemann surfaces. Comm. Pure Appl. Math. 55 (2002), no. 6, 728-771}


\bibitem{12}{B. Gidas, W-Y. Ni, L. Nirenberg. Symmetry and related properties via the maximum principle.  Comm. Math. Phys.  68  (1979), no. 3, 209-243.}

\bibitem{13}{Hofmann, S. Mitrea, M. Taylor, M. Geometric and transformational properties of Lipschitz domains, Semmes-Kenig-Toro domains, and other classes of finite perimeter domains. J. Geom. Anal. 17 (2007), no. 4, 593?647.}

\bibitem{13}{J.M. Lee, T.H. Parker. The Yamabe problem. Bull.Amer.Math.Soc (N.S) 17 (1987), no.1, 37-91.}

\bibitem{14}{YY. Li. Prescribing scalar curvature on $ {\mathbb S}_n $ and related
Problems. C.R. Acad. Sci. Paris 317 (1993) 159-164. Part
I: J. Differ. Equations 120 (1995) 319-410. Part II: Existence and
compactness. Comm. Pure Appl.Math.49 (1996) 541-597.
}


\bibitem{15}{YY. Li. Harnack Type Inequality: the Method of Moving Planes. Commun. Math. Phys. 200,421-444 (1999).}

\bibitem{16}{YY. Li, I. Shafrir. Blow-up Analysis for Solutions of $
-\Delta u = V e^u $ in Dimension Two. Indiana. Math. J. Vol 3, no 4.
(1994). 1255-1270.}

\bibitem{17}{YY. Li, L. Zhang. A Harnack type inequality for the Yamabe equation in low dimensions.  Calc. Var. Partial Differential Equations  20  (2004),  no. 2, 133-151.}

\bibitem{18}{YY.Li, M. Zhu. Yamabe Type Equations On Three Dimensional Riemannian Manifolds. Commun.Contem.Mathematics, vol 1. No.1 (1999) 1-50.}

\bibitem{19}{L. Ma, J-C. Wei. Convergence for a Liouville equation. Comment. Math. Helv. 76 (2001) 506-514.}

\bibitem{20}{I. Shafrir. A sup+inf inequality for the equation $ -\Delta u=Ve^u $. C. R. Acad.Sci. Paris S\'er. I Math. 315 (1992), no. 2, 159-164.}

\bibitem{21}{Y-T. Siu. The existence of Kahler-Einstein metrics on manifolds with positive anticanonical line bundle and a suitable finite symmetry group.  Ann. of Math. (2)  127  (1988),  no. 3, 585-627}

\bibitem{22}{G. Tian. A Harnack type inequality for certain complex Monge-Amp�re equations.  J. Differential Geom.  29  (1989),  no. 3, 481-488.}

\bibitem{22}{L. Zhang. Blowup solutions of some nonlinear elliptic equations involving exponential nonlinearities. Comm. Math. Phys. 268 (2006), no. 1, 105-133.}

\end{thebibliography}
\end{document}